\documentclass[11pt,a4paper,reqno]{amsart}
\usepackage{amsmath}
\usepackage{amsfonts}
\usepackage{amssymb}
\usepackage[utf8]{inputenc}
\usepackage[T1]{fontenc}
\usepackage{amssymb,amsmath,amsthm}
\usepackage{graphicx}
\usepackage{enumerate}
\usepackage{longtable}
\date{}

\DeclareMathOperator{\conv}{conv}
\DeclareMathOperator{\bd}{bd}

\DeclareMathOperator{\dist}{dist}
\DeclareMathOperator{\vol}{vol}
\newtheorem{theorem}{Theorem}[section]
\newtheorem{lemma}[theorem]{Lemma}

\theoremstyle{remark}
\newtheorem{remark}[theorem]{Remark}

\begin{document}
\pagestyle{plain}
\title{A uniform lower bound on the norms of hyperplane projections of spherical polytopes}

\author{Tomasz Kobos}

\address{Faculty of Mathematics and Computer Science \\ Jagiellonian University \\ Lojasiewicza 6, 30-348 Krakow, Poland}

\email{Tomasz.Kobos@uj.edu.pl}


\begin{abstract}

Let $K$ be a centrally symmetric spherical and simplicial polytope, whose vertices form a $\frac{1}{4n}-$net in the unit sphere in $\mathbb{R}^n$. We prove a uniform lower bound on the norms of all hyperplane projections $P: X \to X$, where $X$ is the $n$-dimensional normed space with the unit ball $K$. The estimate is given in terms of the determinant function of vertices and faces of $K$. In particular, if $N \geq n^{4n}$ and $K = \conv \{ \pm x_1, \pm x_2, \ldots, \pm x_N \}$, where $x_1, x_2, \ldots, x_N$ are independent random points distributed uniformly in the unit sphere, then every hyperplane projection $P: X \to X$ satisfies an inequality $\|P\|_X \geq 1+c_nN^{-(2n^2+4n+6)}$ (for some explicit constant $c_n$), with the probability at least $1 - \frac{3}{N}.$ 

\keywords{Minimal projection, Spherical Polytope, Random Polytope, Determinant}
\end{abstract}

\maketitle

%

\section{Introduction}

Let $X$ be a real normed space of dimension at least $2$. A linear and continuous mapping $P: X \to X$ is called a \emph{projection}, if it satisfies the equation $P^2=P$. By a \emph{hyperplane projection} we shall mean a projection with the image of codimension $1$.

Projection is a very old concept in mathematics and a basic notion of the approximation theory, as it provides an approximation of the identity operator on a subspace, by a linear operator defined on the whole space. For this reason, one often seeks for a projection with the smallest possible operator norm, as the smaller norm yields a better approximation. Such a projection $P$ is called a \emph{norm-minimal projection} from $X$ onto the image of $P$.

Norm-minimal projection were studied by a lot of different authors in a great variety of contexts (see for example \cite{baronti}, \cite{basso}, \cite{blatter}, \cite{lewickichalm}, \cite{lewickichalm2} \cite{foucartskrzypek}, \cite{grunbaum}, \cite{kadecsnobar}, \cite{konig}, \cite{konig2}, \cite{konig3}). The so-called \emph{projection constants} were studied most extensively. Projection constants play a profound role in the functional analysis and the local theory of Banach spaces, as they are deeply connected with some other important numerical invariants of Banach spaces. We refer to Chapter 8 in a monograph \cite{tomczak} for a broader picture on the theory of the projection constants. In terms of the norm, the best possible situation happens, when there exists a projection of norm $1$ onto a given subspace. In this case, we say that a subspace is \emph{$1$-complemented} in the given space. 

By the Hahn-Banach theorem, every $1$-dimensional subspace is $1$-complemented. For this reason, we shall call a projection \emph{non-trivial}, if its image has dimension at least $2$, and it is different from the whole space. In a Hilbert space, every subspace is $1$-complemented by means of the orthogonal projection. Conversely, it is well-known, that if every subspace of a given Banach space $X$, satisfying $\dim X \geq 3$, is $1$-complemented, then $X$ is isometric to a Hilbert space (see for example \cite{kakutani}). Still, most of the classical spaces posses some $1$-complemented subspace of dimension at least $2$, even if they are not necessarily Hilbert. For example, in the classical $\ell_p^n$ space, a hyperplane $Y$ is $1$-complemented if and only if $Y$ is the kernel of a functional, which represented as a vector in $\mathbb{R}^n$, has at most two coordinates that are different from zero (see \cite{baronti}). Study of $1$-complemented subspaces of Banach spaces has a long history and there is a large volume of published research on this topic (see \cite{randsurvey}).

Bosznay and Garay proved in \cite{bosznay} that, in the context of a normed spaces of given dimension $n \geq 3$, this is, in fact, a very rare instance, to posses some non-trivial projection with the norm $1$. It turns out, that the set of $n$-dimensional normed spaces, for which every non-trivial projection $P: X \to X$ has norm strictly larger than $1$, is open and dense in the set of all $n$-dimensional normed spaces. This somewhat reminds of the well-known fact, that the set of continuous and nowhere differentiable functions forms an open and dense subset of the set of continuous functions. Moreover, this naturally raises a question of establishing some explicit, uniform lower bound on the norms of projections of a given space, which is strictly greater than $1$. Thinking more globally, it is natural to define the constant $\rho_n$ as the largest positive number, for which there exists an $n$-dimensional normed space $X$, such that every non-trivial projection $P: X \to X$ satisfies $\|P\|_X \geq 1 + \rho_n$. The fact, that $\rho_n$ is positive, follows immediately from the result of Bosznay and Garay and a standard compactness argument. By a result of \cite{lewickichalm}, on every two-dimensional subspace there is always a projection with the norm at most $\frac{4}{3}$, so obviously we have $\rho_n \leq \frac{1}{3}$ for every $n \geq 3$.

It seems that no positive lower bounds on $\rho_n$ are known. This may be related to the fact, that lower bounds for the norms of projections were studied mostly in the case of specific subspaces, rather than uniformly. Nevertheless, some remarkable results related to the uniform lower bounds were obtained. Gluskin in \cite{gluskin2} and Szarek in \cite{szarek} used norms generated by random polytopes to establish such lower bounds, but only for projections with the rank in a specific range. Later, a similar construction was provided also in \cite{latala}. All of these results give estimates of the following type: there exists an $n$-dimensional normed space $X$, such that for every projection $P: X \to X$ with the rank $m$ in an interval of the form $[\alpha n, \beta n]$ (where $0 < \alpha < \beta < 1$ are constants), we have $\|P\|_X \geq C \sqrt{m}$ (for some constant $C$ depending on $\alpha$ and $\beta$). Asymptotically speaking, this is best possible up to a constant, as the famous result of Kadec and Snobar (see \cite{kadecsnobar}) yields the inequality $\|P\|_X < \sqrt{m}$.

A deeply profound role, that random polytopes play in the modern high-dimensional geometry, has been started with a pioneering previous work of Gluskin in \cite{gluskin}, who used them to prove that the asymptotic order of the diameter of the Banach-Mazur compactum is linear. After that, many different important applications of the random objects in the high-dimensional geometry have been established, including the examples above.

It does not seem possible to apply those methods directly to projections with the rank not in the interval of the form $[\alpha n, \beta n]$. In this case, the examples are generally lacking. However, some results were obtained in \cite{kobos} for the case of hyperplane projections. For each $n \geq 3$ let us define a constant $\rho_n^H$ as the largest positive number, for which there exists an $n$-dimensional normed space $X$, such that every hyperplane projection $P: X \to X$ satisfies $\|P\|_X \geq 1 + \rho^H_n$. Obviously, we have $\rho_n^H \geq \rho_n$ for every $n \geq 3$. Moreover, by a result of Bohnenblust (see \cite{bohnenblust}), every hyperplane admits a projection of norm at most $2-\frac{2}{n}$, and therefore $\rho^H_n \leq 1 - \frac{1}{n}$ for every $n \geq 3$. In \cite{kobos} a uniform lower bound on the norms of the hyperplane projections was provided in the case of the space $X$ being a rather general subspace of the $\ell_{2p}^m$ space (where $p \geq 2$ is an integer). In consequence, we have
\begin{equation}
\label{stalaosz}
\rho_n^H \geq  \left (8 \left ( n + 3 \right )^{5} \right )^{-30(n+3)^2},
\end{equation}
for every $n \geq 4$. This implies an asymptotic lower bound on $\rho_n^H$ of the form
\begin{equation}
\label{asymptotyka}
\rho_n^H \geq \exp(-C n^2 \log n),
\end{equation}
for some absolute constant $C>0$.

The aim of this paper, is to study uniform lower bounds for the norms of hyperplane projections in the setting of the spherical polytopes (by which we mean the convex hulls of points lying on the unit sphere). Our main result gives such an explicit, uniform lower bound for a broad class of normed spaces, with the unit ball being a symmetric spherical, simplicial polytope whose vertices form a $\frac{1}{4n}-$net in the unit sphere. If $K \subseteq \mathbb{R}^n$ is a convex polytope, then by a term \emph{face} we shall mean only $(n-1)$-dimensional face (facet) of $K$. We will say that face $F$ of $K$ is given by a vector $f \in \mathbb{S}_n$ (or corresponds to $f$), if $f$ is a unit vector perpendicular to the affine hyperplane containing $F$. We recall that a convex polytope is called \emph{simplicial}, if every face is an $(n-1)$-dimensional simplex. By $\| \cdot \|$ we shall always mean the Euclidean norm in $\mathbb{R}^n$. By $\mathbb{S}_n = \{ x \in \mathbb{R}^n \ : \ \|x\|=1\}$ we denote the Euclidean unit sphere in $\mathbb{R}^n$ (not to be confused with the unit sphere in $\mathbb{R}^{n+1}$). For an $n$-dimensional normed space $X$, the norm of $X$ will be denoted by $\| \cdot \|_X$ and the same symbol will be used for the operator norm of a projection $P:X \to X$. Two faces of $K$ are called \emph{non-neighbouring} if their intersection is empty. For a given $\varepsilon>0$, a set $X \subseteq \mathbb{S}_n$ is called an $\varepsilon$\emph{-net} if for every point $p \in \mathbb{S}_n$, there exists $x \in X$ such that $\|x-p\| \leq \varepsilon$. Throughout the paper, we always assume that $n \geq 3$ is an integer. Our main result goes as follows.

\begin{theorem}[General lower bound for the spherical polytopes]
Let $N$ be a positive integer and $\alpha, \beta$ positive real numbers. Suppose that points $x_1, x_2, \ldots, x_N \in \mathbb{S}_n$ satisfy the following conditions:
\begin{enumerate}
\label{twglowne}
\item Vertices of a convex polytope $K=\conv\{ \pm x_1, \pm x_2, \ldots, \pm x_N\}$ form a $\frac{1}{4n}$-net in $\mathbb{S}_n$.
\item $K$ is a simplicial polytope.
\item For any $1 \leq i_1 < i_2 <  \ldots < i_{n} \leq N$ we have $|\det(x_{i_1}, x_{i_2}, \ldots,  x_{i_{n}})| \geq \alpha$.
\item For any $n$ pairwise non-neighbouring faces of $K$, given by the vectors $f_1, f_2, \ldots, f_{n} \in \mathbb{S}_n$, we have $|\det(f_1, f_2, \ldots, f_{n})| \geq \beta$.
\end{enumerate}
Let $X$ be the $n$-dimensional real normed space with the unit ball $K$. Then, every hyperplane projection $P:X \to X$ satisfies 
$$\|P\|_X \geq 1 + C_n\alpha^2\beta,$$
where
$$C_n= \frac{2^{\frac{3}{2}n-2} \cdot n^{n-4}}{5 \sqrt{n-1}}.$$
\end{theorem}

Informally speaking, the determinant function that appears in the estimate above, could be considered as some kind of a measure of the linear independency of the vertices and facets of $K$. The probability that a random symmetric polytope, with vertices drawn uniformly and independently in $\mathbb{S}_n$, is simplicial and the determinant function does not vanish on any subset of vertices or pairwise non-neighbouring faces is equal to $1$. Thus, if we pick appropriately large number of points at random from the unit sphere, the resulting normed space, has all hyperplane projections with the norm greater than $1$ with a high probability. This is clearly expected by the previously mentioned result of Bosznay and Garay. The question of estimating this probability for some explicit, uniform lower bound on the norms of hyperplane projection arises naturally. Points $x_1, x_2, \ldots, x_N \in \mathbb{S}_n$ will be called \emph{random points}, if they are distributed independently and uniformly in $\mathbb{S}_n$. In the next result, we provide a uniform lower bound for the norms of hyperplane projections, which holds with a large probability for a random spherical polytope. The result states, that for the symmetric convex hull of $N \geq n^{4n}$ random points in $\mathbb{S}^n$, the corresponding normed space has all hyperplane projections with norm greater than $ 1 + c_nN^{-(2n^2+4n+6)}$ (for some specific constant $c_n$), with the probability at least $1 - \frac{3}{N}$.

\begin{theorem}[Lower bound for random spherical polytopes]
\label{twrandom}
Let $N \geq n^{4n}$ be a positive integer. Let $x_1, x_2, \ldots, x_{N} \in \mathbb{S}_{n}$ be random points and let $X=(\mathbb{R}^n, \| \cdot \|_X)$ be the $n$-dimensional normed space with the unit ball $B_X = \conv \{ \pm x_1, \pm x_2, \ldots, \pm x_N \}$. Then, the probability that every hyperplane projection $P:X \to X$ satisfies
$$\|P\|_X \geq 1 + c_nN^{-(2n^2+4n+6)},$$
where
$$c_n=\frac{n^{2n^2+7n-11}}{e^{2n^2+7n+3}}, $$
is at least $1-\frac{3}{N}$.
\end{theorem}

We can say that the result above quantifies the original result of Bosznay and Garay (in the hyperplane setting) in two different ways. It gives a uniform lower bound on the norms of projections, but it also estimates the measure of the spherical polytopes with given number of vertices, which satisfy it. In the three-dimensional case, this gives an estimate of $1 + cN^{-36}$ for some explicit constant $c>0$ -- see Remark \ref{3d}. We   note, that even if we work with random polytopes, our methods are different than those from previously mentioned papers \cite{gluskin2}, \cite{latala}, \cite{szarek}. This may stem from the fact that the hyperplane case seems to be rather dissimilar to the case of projections with the rank depending linearly on $n$. In particular, we do not rely on some more advanced variants of the concentration of measure, but we use only basic probabilistic tools, such as Markov's inequality. Let us also remark, that we took some care, to keep all constants appearing in our estimates as explicit as possible, but in some instances they were estimated rather crudely, in order to keep the results clearer.

Since we consider polytopes approximating the unit sphere very well, the corresponding norm is close to the Euclidean norm and we are working rather locally around it. Therefore, as the orthogonal projection has norm $1$, it is not reasonable to expect that our results will yield an optimal lower bound on the constant $\rho_n^H$.  By taking $N=n^{4n}$ in Theorem \ref{twrandom} we get a lower bound
$$\rho_n^H \geq \exp(-C n^3 \log n),$$
which is worse than the lower bound (\ref{asymptotyka}) obtained in \cite{kobos}. We also note that the estimate (\ref{stalaosz}) from \cite{kobos} works only for $n \geq 4$. Thus, our result gives a first non-trivial bound for the three-dimensional constant $\rho_3^H=\rho_3$. See Remark \ref{3d2} for some numerical estimate on $\rho_3$ that can be deduced with our approach. It should be emphasized, that main goal of the paper is to study uniform bounds on the norms of projections in the setting of the spherical polytopes with large number of vertices. The general discussion on $\rho_n$ and $\rho_n^H$ was presented for the sake of motivating this line of research. We consider the derived bounds on $\rho_n^H$, which are likely far from being optimal, only as an interesting by-product of our investigation.

The paper is organized as follows. In Section \ref{dowodogolne} we prove Theorem \ref{twglowne}. In Section \ref{probabil} we establish Theorem \ref{twrandom} by applying Theorem \ref{twglowne} in combination with several auxiliary results concerning random polytopes. The paper is concluded in Section \ref{concluding}, where some further remarks related to our results, are provided.

\section{Proof of the general lower bound for the spherical polytopes}
\label{dowodogolne}

In this Section we prove Theorem \ref{twglowne}. We start with some simple auxiliary results.

\begin{lemma}
\label{siec}
Let $N$ be a positive integer and suppose that a set $\{x_1, x_2, \ldots, x_N\} \subseteq \mathbb{S}_n$ is an $\varepsilon$-net for some $0<\varepsilon<1$. Then, every face of a convex polytope $K = \conv \{ x_1, x_2, \ldots, x_N \}$ has diameter not greater than $2\varepsilon$ and an inclusion $(1-\varepsilon)\mathbb{S}_n \subseteq K$ holds.
\end{lemma}

\emph{Proof}. Let $F$ be any face of $K$ and let $x$ be the center of the spherical cap determined by $F$. Clearly, $x$ is equidistant to the vertices of $F$. Let us denote this distance by $d$. We have $d \leq \varepsilon$. By the triangle inequality, the distance between any two vertices of $F$ is at most $2d \leq 2 \varepsilon$, which shows the first claim. For the second claim, let $h$ be the distance between $x$ and the hyperplane contaning $F$. Clearly $h \leq d \leq \varepsilon$. Moreover, the hyperplane tangent in $x$ to $S$ is parallel to the hyperplane containing $F$. Since $\frac{1}{1-h} \leq \frac{1}{1-\varepsilon}$, it is clear that $S \subseteq \frac{1}{1-\varepsilon}K$ and the proof is complete. \qed

\begin{lemma}
\label{iloczyn}
Let $x_1, x_2, \ldots, x_n \in \mathbb{S}_n$ be such that $|\det(x_1, x_2, \ldots, x_n)| \geq \alpha$ for some $\alpha>0$. Then, for any $v \in \mathbb{S}_n$ we have $\max_{1 \leq i \leq n} |\langle x_i, v \rangle | \geq \frac{\alpha}{n}.$
\end{lemma}
\emph{Proof.} Assume on the contrary, that for for some $v \in \mathbb{S}_n$ and for each $1 \leq i \leq n$ we have that  $\langle x_i, v \rangle = r_i \in (-\frac{\alpha}{n}, \frac{\alpha}{n}).$ Let $r=(r_1, r_2, \ldots, r_n) \in \mathbb{R}^n$. Then $\|r\| < \sqrt{\frac{n\alpha^2}{n^2}}=\frac{\alpha}{\sqrt{n}}$. Hence, from the Cramer's rule and the Hadamard inequality it follows that
$$|v_1| = \left | \frac{\det(r, x_2, \ldots, x_n)}{\det(x_1, x_2, \ldots, x_n)} \right | \leq \frac{\|r\|}{\alpha} < \frac{1}{\sqrt{n}}.$$

In the same way we prove that $|v_i| < \frac{1}{\sqrt{n}}$ for $2 \leq i \leq n$. Hence
$$1 = v_1^2 + v_2^2 + \ldots + v_n^2 < n \cdot \frac{1}{n} = 1,$$
which gives the desired contradiction. \qed

Now we are ready to prove Theorem \ref{twglowne}.

\emph{Proof of Theorem \ref{twglowne}.}

We denote by $\| \cdot \|_X$ the norm of $X$ and by $\|\cdot\|_{X^{*}}$ the dual norm. By the second assumption and Lemma \ref{siec}, every face of $K$ is an $(n-1)$-dimensional simplex of a diameter not greater than $d = \frac{1}{2n}$. Moreover, the inclusion $\mathbb{S}_n \subseteq \frac{4n}{4n-1}K$ holds. Let $Y \subseteq X$ be an arbitrary $(n-1)$-dimensional subspace. Suppose that $P: X \to X$ is a projection with the image $Y$, that satisfies the inequality $\|P\|_X < 1 + C_n\alpha^2 \beta$. Let us take $w \in \ker P$ with $\|w\|=1$ and suppose that $Y = \{x \in \mathbb{R}^n \ : \ \langle x, v \rangle = 0\}$ with $\|v\|=1$. Assume that $Y$ has non-empty intersection with the faces $F_1, \ldots, F_k, -F_1, \ldots, -F_k$ of $K$, given by the vectors $f_1, \ldots, f_k, -f_1, \ldots, -f_k \in \mathbb{S}_n$. Let us call a face $F_i$ a \emph{bad} face, if there does not exist a vector $z \in Y \cap F_i$ such that $\dist(z, \bd F_i) \geq s$ (the boundary $\bd F_i$ is considered in the affine hyperplane containing $F_i$), where 
$$s = \frac{2^{\frac{n}{2}-1}\alpha^2}{n^2 \cdot \sqrt{n-1} \cdot d^{n-1}}.$$ Otherwise, a face $F_i$ is called a \emph{good} face. We shall prove the following claim.

\textbf{Claim 1}. If $F$ is a bad face, then there exists a vertex $a$ of $F$ such that $| \langle a, v \rangle | < \frac{\alpha}{n}.$

Indeed, let  $F = \conv \{ a_1, a_2, \ldots, a_n\}$ be a bad face. A region 
$$F'=\{ x \in F: \ \dist(x, \bd F) \geq s\}$$
is a simplex, positively homothetic to $F$ (soon, we shall see that $F'$ is non-empty). As the hyperplane $Y$ does not intersect $F'$, the simplex $F'$ lies in one of the open half-spaces determined by $Y$. Without loss of generality, let us assume that a vertex $a_1$ of $F$ lies in the opposite (closed) half-space -- as $Y$ has a non-empty intersection with $F$, there are vertices in both closed half-spaces. Let $f_1$ be the face of $F$ not containing $a_1$ (thus $f_1$ has dimension $n-2)$. Then, it is clear that $Y$ intersects the parallelotope
$$P_1 = \{ x \in F: \dist(x, f) \leq s \text{ for every face } f \neq f_1 \text{ of } F  \}.$$

Let $a_1'$ be a vertex of $F'$ corresponding to $a_1$. Then $a_1' \in P_1$ and $\|a_1-x\| \leq \|a_1-a_1'\|$ for some $x \in P_1 \cap Y$. We shall now prove an upper estimate on the distance $\|a_1-a_1'\|$. Let $0<k<1$ be the homothety ratio of $F$ and $F'$ and let $r$ be the inradius of $F$. The homothety center $c$ is the incenter of both $F$ and $F'$. In particular, $kr + s = r$, which gives us an equality $k = \frac{r-s}{r}$. Furthermore,
$$\|a_1-c\|=\|a_1-a_1'\|+\|a_1'-c\|=\|a_1-a_1'\|+k\|a_1-c\|,$$
which yields
\begin{equation}
\label{odleglosc}
\|a_1-a_1'\|=(1-k)\|a_1-c\|=\frac{s}{r}\|a_1-c\|\leq \frac{ds}{r}.
\end{equation}
We shall now establish a lower bound on the inradius $r$. Let us observe that
$$\vol (\conv \{0, a_1, \ldots, a_n\}) = \frac{|\det(a_1, \ldots, a_n)|}{n!} \geq \frac{\alpha}{n!},$$
but on the other hand,
$$\vol (\conv \{0, a_1, \ldots, a_n\}) = \frac{h}{n} \cdot \frac{r}{n-1}(S_1 + S_2 + \ldots + S_n),$$
where $h$ denotes the distance of the hyperplane determined by $F$ to the origin, and $S_i$ for $1 \leq i \leq n$ denote the $(n-2)$-dimensional volumes of the faces of $F$. Each face of $F$ is an $(n-2)$-dimensional simplex with the edge length not greater than $d$. Thus, it is well-known, that under these constraints, for each $1 \leq i \leq n$, volume $S_i$ is not greater than that of an $(n-2)$-dimensional regular simplex of edge length $d$, which is equal to 
$$\frac{\sqrt{n-1}}{(n-2)! \cdot 2^{\frac{n}{2}-1}} \cdot d^{n-2}.$$
Combining the two previous estimates with an obvious inequality $h < 1$, we get a lower bound
$$r > \frac{(n-1)\alpha}{n!} \cdot \frac{(n-2)! \cdot 2^{\frac{n}{2}-1}}{\sqrt{n-1} \cdot d^{n-2}} = \frac{2^{\frac{n}{2}-1} \alpha}{n \cdot \sqrt{n-1} \cdot d^{n-2}}.$$
Let us note here, that this shows in particular, that the region $F'$ is non-empty, as by the above inequality we clearly have $s < r$. Now, coming back to (\ref{odleglosc}) we obtain
$$\|a_1-a_1'\| \leq \frac{ds}{r} < \frac{ds \cdot n \cdot \sqrt{n-1} \cdot d^{n-2}}{2^{\frac{n}{2}-1}\alpha} = \frac{\alpha}{n}.$$ 
Hence
$$|\langle a_1, v \rangle |=| \langle a_1-x, v \rangle | \leq \|a_1 - x\| \leq \|a_1-a_1'\| < \frac{\alpha}{n},$$
which proves Claim 1.

Now, we shall establish a similar claim for the good faces.

\textbf{Claim 2.} If $F$ is a good face, given by the vector $f \in \mathbb{S}_n$, then $|\langle w, f \rangle | < \frac{\beta}{n}$.

Since $F$ is a good face, there exists a vector $z \in Y \cap F$, such that $\dist(z, \bd F) \geq s$. This means, that a ball $B(z, s)$ with the center $z$ and radius $s$, intersected with the boundary of $K$, is an $(n-1)$-dimensional Euclidean ball contained in $F$. Let us take any real number $\lambda$ satisfying $|\lambda| \leq \frac{s}{5}$. Then, we have
$$\|\|z+\lambda w\|_X z - (z+\lambda w)\|_X = \| \left ( \|z + \lambda w\|_X - 1 \right )z - \lambda w\|_X \leq \left | \|z + \lambda w\|_X - 1 \right | \|z\|_X + |\lambda| \|w\|_X$$
$$= \left | \|z + \lambda w\|_X - \|z\|_X \right |  + |\lambda| \|w\|_X \leq |\lambda| \|w\|_X + |\lambda| \|w\|_X =  2 |\lambda| \|w\|_X  \leq \frac{2 |\lambda|}{1-\frac{1}{4n}} \leq 3|\lambda|.$$
Moreover
$$\|z+\lambda w\|_X \geq 1 - |\lambda| \|w\|_X \geq 1 - \frac{1}{1-\frac{1}{4n}} |\lambda| \geq 1 - 2 |\lambda|.$$

By combining two previous estimates we obtain
$$\left \| z - \frac{z+\lambda w}{\|z+\lambda w\|_X}  \right \|_X \leq \frac{3 |\lambda|}{1-2|\lambda|} \leq s,$$
where the last inequality follows easily from the inequality $|\lambda| \leq \frac{s}{5}$. Finally, we have that
$$\left \|z - \frac{z+\lambda w}{\|z+\lambda w\|_X}\right  \| \leq \left \|z - \frac{z+\lambda w}{\|z+\lambda w\|_X} \right \|_X \leq s.$$
This shows that for $|\lambda| \leq \frac{s}{5}$, the vector $\frac{z+\lambda w}{\|z+\lambda w\|_X}$ belongs to the intersection of a ball $B(z, s)$ with the boundary of $K$ and therefore to the face $F$. In consequence, we have
$$\|z + \lambda w\|_X = \langle z + \lambda w, \tilde{f} \rangle,$$
where $\tilde{f} = \frac{f}{\|f\|_{X^{*}}}$. Thus
$$1 + \lambda \langle w, \tilde{f} \rangle = \langle z + \lambda w, \tilde{f} \rangle = \|z + \lambda w\|_X > \frac{\|P(z + \lambda w)\|_X}{1+C_n\alpha^2 \beta} = \frac{\|z\|_X}{1+C_n\alpha^2 \beta} = \frac{1}{1+C_n\alpha^2 \beta}.$$
Hence,
$$\lambda \langle w, \tilde{f} \rangle \geq \frac{-C_n \alpha^2 \beta}{1 + C_n\alpha^2 \beta} \geq -C_n \alpha^2 \beta$$
By taking $\lambda = \pm \frac{s}{5}$ and using the fact that $\|\tilde{f}\| \geq \|f\|$, we get
$$|\langle w, f \rangle | \leq |\langle w, \tilde{f} \rangle | \leq \frac{5C_n \alpha^2 \beta}{s} = \frac{\beta}{n},$$
by the definitions of $s$ and $C_n$. This proves Claim 2.

Now, with both Claims in our disposal, we can finish the proof of Theorem \ref{twglowne}. We take any point $x \in Y$ with $\|x\|_X=1$ and any two-dimensional subspace $V \subseteq Y$ such that $x \in V$. Let us consider a two-dimensional curve (a broken path), lying in $V \cap \bd K$, which connects $x$ and $-x$. Clearly, its Euclidean length is greater than $2\|x\| \geq \frac{4n-1}{2n} \|x\|_X=\frac{4n-1}{2n}$. This means, that we can find points $p_1, p_2, \ldots, p_{2n-1}$ on this curve such that, for $i \neq j$ we have
$$\|p_i-p_j\| \geq \frac{4n-1}{2n(2n-1)}>\frac{1}{n}=2d.$$
Every point $p_i$ lies in the boundary of $K$, and thus in some face of $K$. Let $F_i$  be any face of $K$ such that $p_i \in F_i$. Note for $i \neq j$, faces $F_i$ and $F_j$ are non-neighbouring. Indeed, if $u \in F_i \cap F_j$, then
$$2d = \frac{1}{n} < \|p_i - p_j\| \leq \|p_i - u\| + \|u-p_j\| \leq 2d,$$
which is a contradiction. It is clear, that in the set $\{ F_1, F_2, \ldots, F_{2n-1} \}$ there are at least $n$ bad faces or at least $n$ good faces. If there are $n$ good faces, then by Claim 1, we get existence of $n$ vertices $a_1, a_2, \ldots, a_n$ of $K$, such that $| \langle a_i, v \rangle | < \frac{\alpha}{n}$ for any $1 \leq i \leq n$. This is an immediate contradiction with Lemma \ref{iloczyn} and the third condition of Theorem. Similarly, if there exist $n$ good faces among $F_1, F_2, \ldots, F_{2n-1}$, then there are $n$ pairwise non-neighbouring faces of $K$, given by the vectors $f_1, f_2, \ldots, f_n \in \mathbb{S}_n$, such that $|\langle  w, f \rangle | < \frac{\beta}{n}$. Again, this contradicts Lemma \ref{iloczyn} combined with the fourth condition of Theorem. This completes the proof. \qed

\section{Random spherical polytopes}
\label{probabil}

In this Section we derive Theorem \ref{twrandom} from Theorem \ref{twglowne}. In order to do this, we need to establish several probabilistic lemmas. These results are strongly based on the ideas developed by different authors. We shall rephrase or modify them, according to our needs. We start with a straightforward reformulation of the lower bound on the probabilistic measure of the spherical cap, that was given in \cite{boroczky}.

\begin{lemma}
\label{czapka}
Let $x \in \mathbb{S}_n$ and $0 < r < 1$. Then, the probabilistic measure of the spherical cap
$$C(x, r) = \{ y \in \mathbb{S}_{n} \ : \ \|x-y\| \leq r \}$$
is at least $\frac{1}{\sqrt{2 \pi (n-1)}} \left ( \frac{r}{\sqrt{2}} \right ) ^{n-1}$.
\end{lemma}

\emph{Proof.} If $ 0 < \varphi \leq \frac{\pi}{2}$ is such that
$$r^2 = 2 - 2\sqrt{1-\sin^2 \varphi},$$
then by the first part of Corollary 3.2 in \cite{boroczky}, it follows that the measure of $C(x, r)$ is at least
$$\frac{1}{\sqrt{2 \pi (n-1)}} \sin^{n-1} \varphi.$$
However
$$r^2 = 2 - 2\sqrt{1-\sin^2 \varphi} = 2 \frac{\sin^2 \varphi}{1 + \sqrt{1- \sin^2 \varphi}} \leq 2 \sin^2 \varphi$$
and hence the claim follows. \qed

In \cite{bourgain} there is an outline of the proof, that for any fixed $\varepsilon>0$, the probability, that $N$ random points on the unit sphere in $\mathbb{R}^n$ form a $N^{-\frac{1}{n-1}+\varepsilon}$-net, tends to $1$, as $N \to \infty$ (see the proof of the $(1.27)$ in Appendix A). In order to obtain more explicit estimate, we give a minor modification of this argument in the following lemma.

\begin{lemma}
\label{probsiec}
Let $N \geq n^{4n}$ be a positive integer. If $x_1, x_2, \ldots, x_{N} \in \mathbb{S}_{n}$ are random points, then the probability, that these random points do not form a $\frac{1}{4n}$-net in the unit sphere, is less than $\frac{1}{N}$.
\end{lemma}

\emph{Proof.} It is well-known, that for any $\varepsilon>0$, there exists a $\varepsilon$-net in the unit sphere of cardinality at least $\left (1 + \frac{2}{\varepsilon} \right )^n$. In particular, there exists a $\frac{1}{8n}$-net of the cardinality $(16n+1)^n \leq (17n)^n$. Hence, let $z_1, z_2, \ldots, z_{(17n)^n}$ be some fixed $\frac{1}{8n}$-net in the unit sphere. For a fixed $1 \leq j \leq (17n)^n$, the probability that each point $x_1, x_2, \ldots, x_N$ is outside the cap $C(z_j, \frac{1}{8n})$, is by Lemma \ref{czapka} at most
$$ \left ( 1 - \frac{1}{\sqrt{2 \pi (n-1)}} \left ( \frac{1}{8\sqrt{2}n} \right ) ^{n-1} \right )^N \leq \left ( 1 - \left ( \frac{1}{8\sqrt{2}n} \right ) ^{n} \right )^N.$$
Therefore, the probability that at least one of the caps $C(z_j, \frac{1}{8n})$ is empty, is at most
$$(17n)^n \left ( 1 - \left ( \frac{1}{8\sqrt{2}n} \right ) ^{n} \right )^N \leq e^{(3+\log n)n} \cdot e^{-Na^{-n}} = e^{3n +  n \log n - Na^{-n}},$$
where $a=8 \sqrt{2}n$. Since $N \geq n^{4n}$ and $n \geq 3$ we have that
$$Na^{-n} \geq \sqrt[4]{N}.$$
It is easy to check that for $N \geq n^{4n}$ and $n \geq 3$ the inequality 
$$3n + n \log n + \log N < \sqrt[4]{N}$$
is true. Hence, the probability that at least one of the caps $C(z_j, \frac{1}{8n})$ is empty is less than $\frac{1}{N}$. It remains to  observe, that if each of these caps is non-empty, then the points $x_i$ form a $\frac{1}{4n}$-net in the unit sphere. This concludes the proof. 
\qed

The last piece of information, that we need to apply Theorem \ref{twglowne}, in order to prove Theorem \ref{twrandom}, is the expected value of the determinant of $n$ random points in $\mathbb{S}_n$. More precisely, we estimate the $(-\frac{1}{2})$-moment of the absolute value of the determinant. We use the fact, that the distribution of the random variable $|\det(x_1, x_2, \ldots, x_n)|$, where $x_i \in \mathbb{S}_n$ are random points, is well-known. We will also refer to the following well-known inequalities on the ratio of Gamma functions:
\begin{equation}
\label{wendelin}
\sqrt{x} \leq \frac{\Gamma\left ( x+1 \right )}{\Gamma\left ( x+\frac{1}{2} \right )} \leq \sqrt{x+\frac{1}{2}} \quad \text{ for every } x>0. 
\end{equation}
The lower bound appeared in this form in \cite{gautschi} and upper bound in \cite{wendel}.

\begin{lemma}
\label{probdet}
Let $M_n$ be the expected value of $|\det(x_1, x_2, \ldots, x_n)|^{-\frac{1}{2}}$, where $x_i\in \mathbb{S}_n$ for $1 \leq i \leq n$ are random points. Then, we have $M_n < e^{\frac{n}{4}}(n-1)$.
\end{lemma}
\emph{Proof.} Let us recall, that a random variable has a \emph{Beta distribution} with parameters $\alpha_1, \alpha_2>0$, if its density is given by
$$g(t)=\frac{\Gamma(\alpha_1+\alpha_2)}{\Gamma(\alpha_1)\Gamma(\alpha_2)}t^{\alpha_1-1}(1-t)^{\alpha_2-1},$$
for $t \in (0, 1)$. It turns out that the random variable $\det(x_1, x_2, \ldots, x_n)^2$ has the distribution $\prod_{i=1}^{n-1} \beta_{\frac{i}{2}, \frac{n-i}{2}}$, where $\beta_{\alpha_1, \alpha_2}$ has a Beta distribution with parameters $\alpha_1, \alpha_2$ and the variables in the product are independent (see \cite{mathai} and \cite{kabluchkodet}). The $(-\frac{1}{4})$-moment of the Beta variable with parameters $\alpha_1>\frac{1}{4}, \alpha_2>0$ is equal to
$$\frac{\Gamma\left ( \alpha_1-\frac{1}{4} \right )\Gamma(\alpha_1+\alpha_2)}{\Gamma(\alpha_1)\Gamma(\alpha_1+\alpha_2-\frac{1}{4})}.$$
Therefore,
$$M_n= \prod_{i=1}^{n-1} \frac{\Gamma \left ( \frac{i}{2}-\frac{1}{4} \right ) \Gamma\left ( \frac{n}{2} \right )}{\Gamma\left ( \frac{i}{2} \right)\Gamma\left( \frac{n}{2} - \frac{1}{4} \right )}$$

From the inequality (\ref{wendelin}) we have
$$\frac{\Gamma\left ( \frac{n}{2} \right)}{\Gamma\left ( \frac{n}{2}-\frac{1}{4} \right )} \leq \left ( \frac{n}{2}-\frac{1}{4} \right )^{\frac{1}{4}} \leq \left (\frac{n}{2} \right )^{\frac{1}{4}}. $$

Similarly,
$$\frac{\Gamma\left ( \frac{i}{2}-\frac{1}{4} \right )}{\Gamma\left ( \frac{i}{2} \right)} \leq \frac{\left ( \frac{i}{2}\right )^{\frac{3}{4}}}{\frac{i}{2} - \frac{1}{4}} \leq \frac{ \left ( \frac{i}{2} \right )^{\frac{3}{4}}}{\left ( \frac{i-1}{2} \right )^{\frac{3}{4}}},$$
for $2 \leq i \leq n-1.$ For $i=1$, we can check by a direct calculation that
$$\frac{\Gamma\left ( \frac{1}{4} \right )}{\Gamma\left ( \frac{1}{2} \right)} \leq \frac{2}{\left ( \frac{1}{2} \right )^{\frac{1}{4}}}.$$
Thus,
$$\prod_{i=1}^{n-1} \frac{\Gamma \left ( \frac{i}{2}-\frac{1}{4} \right )}{\Gamma\left ( \frac{i}{2} \right)} \leq 2 \cdot \left ( \frac{(n-2)!}{2^{n-2}} \right )^{-\frac{1}{4}} \cdot \left ( \frac{n-1}{2} \right )^{\frac{3}{4}} = 2 \cdot \left ( \frac{(n-1)!}{2^{n-1}} \right )^{-\frac{1}{4}} \cdot \left ( \frac{n-1}{2} \right )$$
$$=\left ( \frac{(n-1)!}{2^{n-1}} \right )^{-\frac{1}{4}} \cdot (n-1).$$
Hence, we conclude that
$$M_n \leq \left ( \frac{n^n}{n!} \right )^{\frac{1}{4}} \cdot (n-1) < e^{\frac{n}{4}} \cdot (n-1),$$
by the Stirling's approximation formula. This finishes the proof. \qed

Finally, we are ready to move to the proof of Theorem \ref{twrandom}. We will use the following simple upper bound on the binomial coefficient: $\binom{N}{n} \leq \left ( \frac{Ne}{n} \right )^n$ for $N \geq n$.

\emph{Proof of Theorem \ref{twrandom}.}

Let us consider following probability events:
\begin{enumerate}[(i)]

\item Points $x_1, x_2, \ldots, x_N$ do not form a $\frac{1}{4n}$-net in $\mathbb{S}_n$.

\item $|\det(x_{i_1}, x_{i_2}, \ldots, x_{i_n})| \leq \displaystyle \left ( \frac{e^{\frac{5n}{2}}}{n^{2n-2}}N^{2n+2} \right )^{-1}$ for some indices $1 \leq i_1 < i_2 < \ldots < i_n \leq N$.

\item There exist pairwise non-neighbouring faces $F_1, F_2, \ldots, F_n$ of $B_X$, perpendicular to vectors $f_i \in \mathbb{S}_n$ (where $1 \leq i \leq n$), such that $|\det(f_1, f_2, \ldots, f_n)| \leq \displaystyle  \left ( \frac{e^{2n^2+\frac{5n}{2}}}{n^{2n^2+2n-2}}N^{2n^2+2} \right )^{-1}.$

\end{enumerate}
We shall prove that each of these three events has probability at most $\frac{1}{N}$. We can also assume that $B_X$ is a simplicial polytope, as this happens with the probability $1$.

\begin{enumerate}[(i)]

\item Our claim follows directly from Lemma \ref{probsiec}. 

\item We use Lemma \ref{probdet} and the Markov's inequality applied for the random variable 
$$|\det(y_1, y_2, \ldots, y_n)|^{-\frac{1}{2}},$$
where $y_1, y_2, \ldots, y_n \in \mathbb{S}_n$ are random points. For any $a>0$ we have:
$$P \left ( |\det(y_1, y_2, \ldots, y_n)|^{-\frac{1}{2}} \geq a (n-1) e^{\frac{n}{4}} \right ) \leq a^{-1},$$
which is equivalent to
\begin{equation}
\label{prawd}
P \left ( |\det(y_1, y_2, \ldots, y_n)| \leq \left ( a (n-1) e^{\frac{n}{4}} \right )^{-2} \right ) \leq a^{-1}.
\end{equation}
Using this inequality for $a=N \binom{N}{n}$, with the union bound for all possible choices of $n$ points from $x_1, x_2, \ldots, x_N$, along with an estimate $\binom{N}{n} \leq \left ( \frac{Ne}{n} \right )^n$, we get the desired upper bound of $\frac{1}{N}$.

\item If $A \subseteq \{1, 2, \ldots, N\}$ is an $n$-element set, then by $x_A^{\perp} \in \mathbb{S}_n$ we denote a unit vector perpendicular to the affine hyperplane containing the points $x_i$ for $i \in A$ (this hyperplane is determined uniquely with the probability $1$). To be more precise, we can define $x^{\perp}_A$ as the exterior product of $x_i$ for $i \in A$. A set of vectors $\{f_1, f_2, \ldots, f_n\} \subseteq \mathbb{S}_n$, corresponding to a set of $n$ pairwise non-neighbouring faces of $B_X$, can be represented as $\{x_{A_1}^{\perp}, x_{A_2}^{\perp}, \ldots, x_{A_n}^{\perp} \}$ for a certain $n$-element family $\{A_1, A_2, \ldots, A_n\}$ of $n$-element disjoint subsets of the set $\{1, 2, \ldots, N\}$. Let us denote by $\mathcal{F}$ the collection of all such families. Then, for a given $c>0$, the probability that there exist pairwise non-neighbouring faces of $B_X$, given by vectors $f_1, f_2, \ldots, f_n \in \mathbb{S}_n$, that satisfy $|\det(f_1, f_2, \ldots, f_n)| \leq c$ is, by the union bound, at most
$$P\left ( \bigcup_{\{A_1, \ldots, A_n\} \in \mathcal{F}} |\det(x_{A_1}^{\perp}, x_{A_2}^{\perp}, \ldots, x_{A_n}^{\perp})| \leq c) \right )$$
$$\leq\sum_{\{A_1, \ldots, A_n\} \in \mathcal{F}} P \left ( |\det(x_{A_1}^{\perp}, x_{A_2}^{\perp}, \ldots, x_{A_n}^{\perp})| \leq c) \right )$$
$$=\# \mathcal{F} \cdot  P \left ( |\det(x_{B_1}^{\perp}, x_{B_2}^{\perp}, \ldots, x_{B_n}^{\perp})| \leq c) \right ),$$
where $\{B_1, B_2, \ldots, B_n\} \in \mathcal{F}$ is any fixed family (this probability is clearly the same for any family in $\mathcal{F}$). Since the points $x_i$ are chosen uniformly and independently and the sets $B_i$ are disjoint, it easily follows that $\{x_{B_1}^{\perp}, x_{B_2}^{\perp}, \ldots, x_{B_n}^{\perp}\}$ is a set of $n$ independent random points in $\mathbb{S}_n$, with respect to the uniform distribution.  Thus, we can use the Markov's as in (\ref{prawd}) for $a = \# \mathcal{F} \cdot N$ to get
$$P \left ( |\det(x_{B_1}^{\perp}, x_{B_2}^{\perp}, \ldots, x_{B_n}^{\perp})| \leq \left ( \#\mathcal{F} \cdot N (n-1) e^{\frac{n}{4}} \right )^{-2} \right ) \leq (\# \mathcal{F} \cdot N)^{-1}.$$
Using a bound $\binom{N}{n} \leq  \left ( \frac{Ne}{n} \right )^n$ and the Stirling's approximation formula, we can estimate the cardinality of $\mathcal{F}$ as follows
$$\# \mathcal{F} = \frac{1}{n!} \cdot \binom{N}{n} \cdot \binom{N-n}{n} \cdot \ldots \cdot \binom{N-(n-1)n}{n} \leq \left ( \frac{e}{n} \right )^n \binom{N}{n}^n$$
$$\leq \left ( \frac{e}{n} \right )^n \cdot \left ( \frac{Ne}{n} \right )^{n^2} = N^{n^2} \left ( \frac{e}{n} \right )^{n^2+n}.$$
Hence
$$P \left ( |\det(x_{B_1}^{\perp}, x_{B_2}^{\perp}, \ldots, x_{B_n}^{\perp})| \leq \left (  \left ( \frac{e^{n^2+\frac{5n}{4}}}{n^{n^2+n-1}} \right) N^{n^2+1} \right )^{-2} \right ) \leq (\# \mathcal{F} \cdot N)^{-1}$$
and in consequence, the probability that  there exist pairwise non-neighbouring faces of $B_X$ $\{f_1, f_2, \ldots, f_n\} \subseteq \mathbb{S}_n$ satisfying 
$$|\det(f_1, f_2, \ldots, f_n)| \leq \left ( \frac{e^{2n^2+\frac{5n}{2}}}{n^{2n^2+2n-2}}N^{2n^2+2} \right )^{-1}$$
is at most $\frac{1}{N}$.
\end{enumerate}

The result follows now from Theorem \ref{twglowne} and easy computation, where 
$$\alpha=  \left ( \frac{e^{\frac{5n}{2}}}{n^{2n-2}}N^{2n+2} \right )^{-1} \: \text { and } \: \beta= \left ( \frac{e^{2n^2+\frac{5n}{2}}}{n^{2n^2+2n-2}}N^{2n^2+2} \right )^{-1}.$$ \qed

\section{Concluding remarks}
\label{concluding}

In the following section we present some remarks related to previous results. We start with the three-dimensional setting.

\begin{remark}
\label{3d}
In the three-dimensional case, the uniform estimate given in Theorem \ref{twrandom} gives us
$$\|P\|_X \geq 1 + cN^{-36},$$
where $c= \left ( \frac{9}{e^3} \right )^{14}$. By taking $N=3^{12}$ we get a first non-trivial lower bound on $\rho_3=\rho^H_3$. Better numerical bound on $\rho_3$ will be given in the next remark.
\end{remark}

\begin{remark}
\label{3d2}
It is not hard to prove, that in the three-dimensional case, the condition $(1)$ in Theorem \ref{twglowne} can be replaced with an assumption, that the volume of $K$ is greater than $4$ (here we mean the standard volume in $\mathbb{R}^3$) and length of every edge of $K$ is less than $\frac{1}{4}$. With a help of a computer program, a spherical polytope $K$, satisfying these conditions was found. Number of vertices of $K$ is equal to $434$ and $\alpha = 5.303 \cdot 10^{-7}$, $\beta = 1.244 \cdot 10^{-7}$. This gives us a  better numerical estimate than in Remark \ref{3d}:
$$\rho_3 > 932 \cdot 10^{-23}.$$
It is rather hard to believe, that this estimate could be close to the true value of $\rho_3$. Still, in the class of spherical polytopes with $434$ vertices, it is not necessarily so weak.
\end{remark}

It turns out that if the dimension is fixed, the polynomial bound, given in terms of $N$, can be improved.

\begin{remark}
\label{asymptotics}
For the asymptotics with $n$ fixed and $N \to \infty$, the estimate given in Theorem \ref{twrandom} can be strengthened to $1 + c_{n, \varepsilon}N^{-(n^2+3n+3+\varepsilon)}$ for any $\varepsilon>0$ and some constant $c_{n, \varepsilon}$, depending on $n$ and $\varepsilon$. Indeed, the expected value of a random variable $|\det(x_1, x_2, \ldots, x_n)|^{-1+\varepsilon}$ is finite for every $\varepsilon>0$. This can be easily deduced from the distribution of a random variable $\det(x_1, x_2, \ldots, x_n)^2$, given in Lemma \ref{probdet}. Thus, we can replace the exponent of $-2$, in the right hand sides of the inequalities given in the properties $(ii)$ and $(ii)$ in the proof of Theorem \ref{twrandom}, to the exponent of $-1-\varepsilon$, albeit with some different constant depending on $n$ and $\varepsilon$.
\end{remark}

\begin{remark}
\label{l1}
In Theorem \ref{twglowne} it seems to be essential that the vertices of $K$ form a good approximation of the unit sphere (the first condition). Let us consider a convex polytope $K=\conv\{\pm e_1, \pm e_2, \ldots, \pm e_n\}$, where $e_i$ is the $i$-th vector from the canonical unit basis of $\mathbb{R}^n$. Thus, the convex polytope $K$ is the unit ball of the $\ell_1$ norm in $\mathbb{R}^n$. In this case we have:
\begin{itemize}
    \item $K$ is a simplicial polytope, that is, the second condition of Theorem \ref{twglowne} is satisfied.
    \item The third condition is satisfied with $\alpha=1$.
    \item The fourth condition is satisfied with any positive $\beta$, since it is not possible to choose any $n$ pairwise non-neighbouring faces of $K$.
\end{itemize}
Yet, in the $\ell_1$ norm there are hyperplane projections of norm $1$ -- for example any projection $P: \mathbb{R}^n \to \mathbb{R}^n$ of the form
$$P(x) = x - x_i e_i,$$
where $1 \leq i \leq n$. This shows that without the first condition of Theorem \ref{twglowne} (or some substitute of it), it is not possible to give a uniform lower bound on the norms of the hyperplane projections in terms of the minimum of the determinant function of vertices and faces.

This also explains, why is it necessary to assume that the number of random points is large enough in Theorem \ref{twrandom}, that is $N \geq n^{4n}$. By a standard volume argument, yielding a lower bound on the cardinality of an $\varepsilon$-net in the unit sphere, it is not possible to significantly improve the estimate given in Lemma \ref{probsiec}. Therefore, the presented approach works only for spherical polytopes with a sufficiently large number of vertices.

\end{remark}

\begin{remark}
\label{upper}
We were interested only in lower bounds on the norms of the hyperplane projections. In the context of random symmetric spherical polytopes with $2N$ vertices, establishing some good upper bounds on the projection constants of hyperplanes, also seems to be an interesting and non-trivial problem. Here we briefly explain, how to get a simple uniform upper bound that goes to $1$ and that holds with probability tending to $1$, as $N \to \infty$. Indeed, by a result of \cite{bourgain}, we know that $N$ random points in $\mathbb{S}^n$, form an $N^{-\frac{1}{2(n-1)}}-$net with the probability tending to $1$, as $N \to \infty$ (see the comment before Lemma \ref{probsiec}). If this condition is satisfied, then from Lemma \ref{siec} it follows that $\left (1-N^{-\frac{1}{2(n-1)}} \right) \mathbb{S}_n \subseteq K$, where $K = \conv \{\pm x_1, \ldots, \pm x_N \}$. Hence, if we denote by $X=(\mathbb{R}^n, \| \cdot \|_X)$ the normed space, for which $B_X=K$, then
$$\|x\|_2 \leq \|x\|_X \leq \frac{1}{1-N^{-\frac{1}{2(n-1)}}} \|x\|_2 \leq \left ( 1 + 2N^{-\frac{1}{2(n-1)}} \right ) \|x\|_2.$$
Now, let $Y \subseteq X$ be an arbitrary hyperplane. If $P: \mathbb{R}^n \to Y$ is an orthogonal projection, then $\|P\|_2=1$ and from the estimates above it easily follows that $\|P\|_X \leq 1 + 2N^{-\frac{1}{2(n-1)}}$. Thus, combining this with Theorem \ref{twrandom} we get, that for $N \to \infty$, the projection constant of an arbitrary hyperplane in $X$ lies in the interval
$$[1 + c_nN^{-(2n^2+4n+6)}, 1+2N^{-\frac{1}{2(n-1)}}],$$
with the probability tending to $1$. We do not know, to what extent this range is best possible, but it is likely that it could be improved (restricted) significantly.
\end{remark}

We must note, that the asymptotic lower bound (\ref{asymptotyka}) on $\rho_n^H $ does not seem to be optimal. The same is true for the estimate $\rho_3 > 175 \cdot 10^{-16}$ on the three-dimensional constant. Moreover, we do not have any non-trivial bound on the constant $\rho_n$ for $n \geq 4$. 

Providing some different example of a class of normed spaces, for which all non-trivial projections/hyperplane projections satisfy some explicit uniform lower bound on the norm, is of its own interest, even if it does not lead to improvement in the global estimates. We do not know if our results for the random spherical polytopes are optimal. Studying this example in more depth, certainly seems to be interesting. Considering the importance of random constructions in modern functional analysis, it is reasonable to believe, that random polytopes could possibly yield much better bounds on the constants $\rho_n$ and $\rho_n^H$. It is likely, that our methods could be extended to all non-trivial projections, providing some bound on $\rho_n$. However, to obtain a better bound on $\rho_n^H$, it would be probably necessary to use random polytopes with the number of vertices of a smaller order than $n^{4n}$. Even in the hyperplane setting, this would possibly require some completely new ideas.

\section{Acknowledgements}

We are grateful to Zakhar Kabluchko for an extensive discussion concerning the probabilistic part of the paper and to Simon Foucart for useful suggestions. We also thank Marin Varivoda for a help with a computer assisted computation, that resulted in a numerical bound presented in Remark \ref{3d2}.


\begin{thebibliography}{}
\bibitem{baronti} Baronti, M., Papini, P.L.: Norm-one projections onto subspaces of $\ell_p$. Annali di Mat. Pura ed Appl. \textbf{152}, 53-61 (1988)
\bibitem{basso} Basso, G.: Computation of maximal projection constants. J. Funct. Anal. \textbf{277}, 3560-3585 (2019)
\bibitem{blatter} Blatter, J., Cheney, E.W.: Minimal projections on hyperplanes in sequence spaces. Annali di Mat. Pura ed Appl. \textbf{101}, 215-227 (1974)
\bibitem{bohnenblust} Bohnenblust, F.: Subspaces of $\ell_p^n$ spaces. Amer. J. Math \textbf{63}, 64-72 (1941)
\bibitem{bosznay} Bosznay, A.P., Garay B.M.: On norms of projections. Acta Sci. Math. \textbf{50}, 87-92 (1986)
\bibitem{bourgain} Bourgain J., Sarnak P., Rudnick Z.: Local Statistics of Lattice Points on the Sphere. Modern Trends in Constructive Function Theory. Contemp. Math. \textbf{661}, 269-282 (2016)
\bibitem{boroczky} B\"or\"oczky K., Wintsche G.: Covering the Sphere by Equal Spherical Balls. In: Discrete Comput. Geometry- The Goldman-Pollack Festschrift \textbf{25}, B. Aronov, S. Baz\'u, M. Sharir, J. Pach (eds.) pp. 237-253. Springer (2003) 
\bibitem{lewickichalm} Chalmers B.L., Lewicki G.: A proof of the Gr\"unbaum conjecture. Studia Math. \textbf{200}, 103-129 (2010) 
\bibitem{lewickichalm2} Chalmers B.L., Lewicki G.: Three-dimensional subspace of $l_{\infty}^5$ with maximal projection constant. J. Funct. Anal. \textbf{257}, 553-592 (2009)
\bibitem{foucartskrzypek} Foucart S., Skrzypek L.: On maximal relative projection constants. J. Math. Anal. Appl. \textbf{447}, 309-328 (2017)
\bibitem{gautschi} Gautschi W.: Some Elementary Inequalities Relating to the Gamma and Incomplete Gamma Function. J. Math. and Phys. \textbf{38}, 77-81 (1959)
\bibitem{gluskin} Gluskin E.D.: The diameter of the Minkowski compactum is approximately equal to $n$. Funct. Anal. Appl. \textbf{15} 72-73 (1981) (In Russian)
\bibitem{gluskin2} Gluskin E.D.: Finite-dimensional analogues of spaces without a basis. Dokl. Akad. Nauk SSSR \textbf{261}, 1046–1050 (1981) (In Russian)
\bibitem{grunbaum} Gr\"unbaum B.: Projection constants. Trans. Amer. Math. Soc. \textbf{7}, 451-465 (1960)
\bibitem{kabluchkodet} Grote J., Kabluchko Z., Th\"ale C.: Limit theorems for random simplices in high dimensions, ALEA, Lat.  Am.  J. Probab.  Math.  Stat. \textbf{16}, 141-177 (2019)
\bibitem{kadecsnobar} Kadec M.I., Snobar M.G.: Certain functionals on the Minkowski compactum. Mat. Zametki \textbf{10}, 453-457 (1971)
\bibitem{kobos} Kobos T.: A uniform estimate of the relative projection constant. J. Approx. Theory \textbf{225}, 58-75 (2018)
\bibitem{kakutani} Kakutani S.: Some characterizations of Euclidean space. Jpn. J. Math. \textbf{16}, 93-97 (1939)
\bibitem{konig} K\"onig H., Lewis D.R., Lin P.K.: Finite-dimensional projection constants. Studia Math. \textbf{75}, 341-358 (1983)
\bibitem{konig2} K\"onig H., Tomczak-Jaegermann N.: Norms of minimal projections. J. Funct. Anal. \textbf{119}, 253-280 (1994)
\bibitem{konig3} K\"onig H.: Spaces with large projection constants. Israel J. Math \textbf{50}, 181-188 (1985)
\bibitem{latala} Latala R., Mankiewicz P., Oleszkiewicz K., Tomczak-Jaegermann N.: Banach-Mazur distances and projections on random subgaussian polytopes. Discrete Comput. Geom. \textbf{38}, 29-50 (2007)
\bibitem{mathai} Mathai A.M.: Distributions of random volumes without using integral geometry techniques. In: Probability and Statistical Models with Applications, Chapman and Hall/CRC (2001)
\bibitem{randsurvey} Randrianantoanina B.: Norm one projections in Banach spaces.Taiwaneese J. Math. \textbf{5}, 35–95 (2001)
\bibitem{szarek} Szarek S.J.: The finite dimensional basis problem, with an appendix on nets of Grassman manifold. Acta Math. \textbf{159}, 153–179 (1983)
\bibitem{tomczak} Tomczak-Jaegermann N.: Banach-Mazur distances and finite-dimensional operator ideals, Pitman Monographs and Surveys in Pure and Applied Mathematics \textbf{38}, New York (1989)
\bibitem{wendel} Wendel J.G.: Note on the gamma function. Amer. Math. Monthly \textbf{55}, 563-564 (1948)

\end{thebibliography}
\end{document}